\documentclass[12pt]{article}
\usepackage{amsfonts}
\usepackage{amssymb}
\usepackage{mathtools}

\newtheorem{theorem}{Theorem}
\newtheorem{theorem*}{Theorem}

\newtheorem{lemma}{Lemma}
\newtheorem{proposition}{Proposition}
\newtheorem{definition}{Definition}
\newtheorem{definition*}{Definition}
\newtheorem{remark}{Remark}
\newtheorem{example}{Example}
\newtheorem{propositiondefinition}{Proposition-Definition}

\begin{document}

\begin{flushright}
{\bf To the memory of O.N. Vvedenskii}
\end{flushright}

\begin{center}
{{\Large \bf  Abelian varieties,  homogeneous spaces and duality. I} \\
(with  Mass Formulas, Formal Groups and Shtukas)
}
\end{center}

\begin{center}
{\bf N.M. Glazunov } \end{center}

\begin{center}
{\rm Glushkov Institute of Cybernetics NASU, Kiev, } \\
{\rm Institute of Mathematics and Informatics Bulgarian Academy of Sciences }\\
{\rm  Email:} {\it glanm@yahoo.com }
\end{center} 

\bigskip

{\bf Abstract.}
{\cal The article is dedicated to the memory of
 Oleg Nikolaevich Vvedenskii (1937-1981).
The results obtained by O. N. Vvedenskii are presented, as well as selected new results
 of the authors, which develop the study of arithmetic algebraic geometry in the directions of
 crystalline cohomology, fundamental groups of schemes, torsors, dualities.
 Elements of ontology of Vvedenskii's research  are also given.
A continuation of the review of Vvedenskii's results, as well as a review of new selected results, including variants of Smith–Minkowski–Siegel mass formula and Drinfeld shtukas, will be presented in the second part of the paper. }

\bigskip


\bigskip

{\bf Thanks.}{ The author is deeply grateful to the Bulgarian Academy of Sciences, the Institute of Mathematics and Informatics of the Bulgarian Academy of Sciences, Professor P. Boyvalenkov for their support.\\
The author was supported by Simons grant 992227.}

\bigskip

``I remember well, for example, stories about a seminar led by A.0. Gelfond,
B.N. Delone and L.G. Shnirelman, where they tried to understand the  class field theory and came to the conclusion that it was hopeless. Later, as a student, I participated in a seminar by B.N. Delone and A.G. Kurosh on the same topic, which ended with the same result. With regard to algebraic geometry, and especially the works of the Italian school, even such attempts were not made - the belief was widespread that it was impossible to understand them.'' ( I. Shafarevich, 1990)

\section{Introduction}

The article is dedicated to the memory of
 Oleg Nikolaevich Vvedenskii (1937-1981).
O.N. Vvedenskii 
 was a student of Academician I.R. Shafarevich.
 Vvedenski's research and his received
 results are connected with duality in elliptic curves and with the corresponding Galois cohomology over local fields, with the Shafarevich-Tate pairing, and
with other pairings, with  local and quasi-local class field  theories of elliptic curves, 
 with the theory of Abelian varieties of dimension greater than 1, with the theory of commutative formal groups over local 
 fields and over  quasi local fields, with M. Artin effect in abelian varieties.
Both results obtained by O. N. Vvedenskii and selected new results of authors are presented that develop the study of arithmetic algebraic geometry in the directions of crystalline cohomology, fundamental groups of schemes, torsors, dualities.
A continuation of the review of Vvedenskii's works, as well as a review of new selected results, including variants of Smith–Minkowski–Siegel mass formula and Drinfeld shtukas, will be presented in the second part of the paper.
Let $S$ be a scheme. An Abelian scheme over $S$ is a $S$-group scheme $A \to S$ which
proper, flat, finitely presented, and which has smooth and connected geometric fibers.

Professor I. Dolgachev is the I. Shafarevich's student. Vvedenskii met with I. Dolgachev. Prof. Dolgachev supported the work of Vvedenskii's students and acted as an opponent in the defense of their dissertations.
The monograph \cite{cdlk} by Dolgachev and co-authors contains references to the works of Vvedenskii. We also note that the monograph \cite{cdlk} gives some answer to the  last phrase
the above quote by I. Shafarevich.

In the monograph by J. Milne \cite{milne1980}, translated into Russian, edited by I.R. Shafarevich on the initiative of O. N. Vvedenskii (who was also one of the translators), along with  the terminology of principal homogeneous spaces, the concept of a torsor, or, more precisely, a $G$-torsor, is used.
Note that the well-known monograph ``Arithmetic Duality Theorems'' by Milne \cite{milne2006} lists the main works of 
O.N. Vvedenskii on the arithmetic duality theory.
Since presently the presentation of the corresponding results in the language of torsors has become generally accepted [15, 16, 19, 20, 21, 22], and we use, along with the term "principal homogeneous space", the language of torsors.
Elements of ontology of Vvedenskii's research  are also given.

A continuation of the review of Vvedenskii's results and result of his students \cite{Gl73,GlFG,Kon1,Kon2,Gl1,Gl2,Gl3,Gl4}, as well as a review of new selected results, including variants of Smith–Minkowski–Siegel mass formula and Drinfeld shtukas, will be presented in the second part of the paper

\section{Elements of Ontology of  O.N. Vvedenskii`s research}

\subsection{List of symbols that O.N. Vvedenskii have used}

$K$   a local field \\
$\overline K$ algebraic closure of $K$ \\
$\mathsf C$   the completion of $\overline K$\\
$U_K$ group of units of $K$ \\
${\mathfrak D}_K$  ring of integers of $K$ \\
$t$ (sometimes $T$) prime element of $K$ \\
$k$  residue field of $K$ \\
$\overline k$ the algebraic closure of $k$ \\
$K{\overline k}$ compositum of fields or composite field  (if exists), \\
$p > 3$ the characteristic of the residue field \\
$G_a$  additive  group of the residue field $k$\\
$G_k$ multiplicative group of the residue field $k$\\
$A$ an Abelian variety  \\
$A_K$ an Abelian variety defined  over $K$\\
$\overline A$  the Picard variety of the variety $A$ \\
$A^{'}$ reduction of  abelian variety  $mod \; t$\\
$A_K^0$ subgroup of points of $A_K$ which reduced to non singular points of $A^{'}$ \\
$\pi_0(A_K) = A_K/A_K^0$ zero dimensional homotopy group \\
$\pi_1(A_K)$ fundamental group of the pro algebraic group of $A_K$ \\
$\Gamma_K$ the kernel of the epimorphism of the reduction $A_K \to A^{'}$\\
$\Gamma_K = \Gamma_K^1 \supset \Gamma_K^2 \supset \ldots, \; \bigcap \Gamma_K^n = 0$ filtration on $\Gamma_K$\\
$L$ a finite normal extension of $K$ \\
${\hat{\mathbb Z}}$ profinite completion of the ring of integers \\
${\mathfrak G}, {\mathfrak g}$ the Galois group of the extension $L/K$\\
$\mathfrak G_K$ the Galois group of the maximal abelian extension of the local field $K$ \\
${\mathfrak G}_{\overline K}$  the Galois group of the algebraic closure of $K$ \\
$H^1(K,A)$  the group of principal homogeneous spaces over $A$, $K$  the quasi local field.\\
$H^0({\mathfrak G},L)$  the zero cohomology group of the group ${\mathfrak G}$ with coefficients in $L$, modified by Tate \\
$H^1({\mathfrak G},L)$    the first cohomology group of the group ${\mathfrak G}$ with coefficients in $L$, modified by Tate\\

{\bf Remark}. In papers by O. Vvedenskii these symbols can have  and another meaning. In the latter case, the meaning of the symbol is specified.\\

\subsection{Concepts and some definitions}

Finite groups and finite group schemes \\
Pontryagin duality \\
Cartier duality \\
Algebraic, quasi algebraic and pro algebraic groups \\
Elliptic curves \\
Neron model \\
Hasse invariant \\
Finite, local,  quasi local, quasi finite, global and quasi global fields \\
Principal homogeneous spaces and abelian varieties \\
Duality in elliptic curves over a local field \\
On the Galois cohomology of elliptic curves defined over a local field \\
Abelian $l$-adic representations and elliptic curves (by Serre) \\
Local class field theory \\
Quasi local class field theory \\
Artin effect in elliptic curves and in abelian varieties. \\
Abelian varieties and formal groups \\
${\mathbb G}_m(R)$ - multiplicative group over a commutative ring $R$. \\
$S$ or $\mathcal S$ - complete defining set of the group $G$ \\
$\lim \; inv \; G/H$ - inverse, or projective limit of groups $G/H \; (H \in S)$.   \\
Etale (étale) topology \\
Etale (étale) sheaf\\

A local field $K$, i.e., a discretely normed field complete with respect
    to the topology induced by the norm. Below we mainly will consoder non archimedean local fields with finite residue field $k$
 and with    normalized discrete valuation   $\nu$  that is the homomorphism  $\nu: K^{*} \to {\mathbb Z}$ is surjective.\\
   Denote by $\#S$ the number of elements of a finite set $S$.    Put $q = p^n = \#k$.
   There are \\  a)  non archimedean local fields $K$ of characteristic $0$: 
   these are finite extensions of $p$-adic fields ${\mathbb Q}_p$; if $[K:{\mathbb Q}_p] = n$ then $n = f\cdot e$ where $f$ is the residue degree and $e$ is the ramification index $\nu(p)$; and \\
   b) the equal characteristic case, when $char \: K = p > 0$ amd $K$ is isomorphic to a field $k((T))$ of formal power series, where $T$     is a uniformization paremeter.\\
    {\bf Remark}. Biside with this definition of the local field O.N. Vvedenskii  (and other reseachers) subdivide non archimedean local fields  on two classes: non archimedean local fields, if the residue field is finite, and  quasi local fields when the residue field is algebraically closed. \\
O.N. Vvedenskii also uses and investigates in the framework of his research  common local fields - a complete discrete valued fields  with a quasi-finite residue field of positive characteristic. Recall that a field $k$ is called quasi-finite if it is perfect and if 
${\mathcal Gal}_k \simeq {\hat {\mathbb Z}}$  where ${\mathcal Gal}_k$ is the
Galois group of the algebraic closure $k_c$ over $k$ and $\hat {\mathbb Z}$ is the completion of the
additive group of the rational integers.\\

The Hasse or Hasse–Witt invariant $H$ is the rank of the Hasse–Witt matrix of a non-singular algebraic curve 
 over a finite field. In the case of elliptic curves it is equal $0, \; (H = 0)$ if the elliptic curve is
 super singular and $H = 1$ if the elliptic curve is ordinary.\\

 According to Weil \cite{weil1}, the principal homogeneous space over $A$ is the algebraic variety $V$ on which $A$ acts
as a group of regular mappings, and the following conditions are satisfied: \\
1) For any $ u, v \in V$ equation
\begin{equation}
au = v 
   \label{e1}
\end{equation}
has a unique solution $a \in A$. \\
2) Mapping that associates a pair $u$ and $v$ solution $a$
equation (\ref{e1}), is a regular mapping of $V \times V$ into $A$.\\

Let $G$  be a group acting on a set $X$. The action is said to be simply transitive if it is transitive and for all $x, y \in X$
 there is a unique $g \in G$  such that $g \cdot x = y$.\\
 
 Let $G$ a smooth algebraic group. 
 A $G$-torsor or a principal $G$-bundle $P$ over a scheme $X$ is a scheme  with an action of $G$ that is locally trivial in the given Grothendieck topology.\\
 
$G$-torsor as  a principal homogeneous space: a $G$-torsor $P$ on a scheme $X$ is a principal homogeneous space for the group scheme
 ${\displaystyle G_{X}=X\times G}$
 (i.e., ${\displaystyle G_{X}}$
acts simply transitively on ${\displaystyle P}$.) \\

\subsubsection{Divisor equivalences} 
Linear equivalence  of divisors  (\cite{Hartshorne}, p.57)  (connected with Picard group of algebraic variety $X$). Two divisors $D$ and $D^{'}$ are said to be linearly equivalent, written $D \sim  D^{'}$, if 
 $D - D^{'}$ is a principal divisor.\\
Algebraic equivalence of divisors (\cite{Hartshorne},  p. 140)    $D \sim_{alg}  D^{'}$,   or $D \equiv D^{'}$ (although $\equiv$ is used more often for
numerical equivalence of divisors). \\
Numerical equivalence of divisors  (\cite{Hartshorne}, p.364)  $\equiv$

\subsection{Problems and conjectures}
Vvedenskii's works \cite{Vvedenskii64,Vvedenskii66,Vvedenskii70,Vvedenskii73a,Vvedenskii74,Vvedenskii76,Vvedenskii79,Vvedenskii79a,Vvedenskii81} are related to the study of the following problems and hypotheses:\\

Duality in elliptic curves over a local field \\

 Galois cohomology of elliptic curves defined over a local field \\
 
 Elliptic class field theory.\\
In a number of works by J. Tate, I. R. Shafarevich, J. Cassels  
and other authors, it turned out that for elliptic curves (Abelian varieties of arbitrary 
dimension) must take place some analog of the classical class field theory 
 of the multiplicative group, the core of which is the duality between the group of principal 
homogeneous spaces over an elliptic curve (an Abelian variety) and some "arithmetic" object 
associated with this curve (variety).\\

 “Universal norms” of formal groups defined over the ring of a local field \\
  
 Duality in elliptic curves over a quasilocal field \\
 
 Pairings in elliptic curves over global fields  \\
 
 The Artin effect in elliptic curves.
 
  Let $A$ be an Abelian variety over a quasi-global field $K$
(i.e. over the field of algebraic functions of one variable with an algebraically closed field of constants $k$). Let 
$char \; K=  char \; k = p>0$. Let $st$
  be the well-known Shafarevich-Tate group corresponding to $A$ considered over $K$. 
Further, for an Abelian group $X$ and a prime number $q$, denote $X_q = Ker(X \stackrel{q}\to X)$.
Similar notation is then retained for the cases
when $X$ is a commutative group scheme or an Abelian sheaf on some Grothendieck topology.
From the works of I. R. Shafarevich, A. P. Ogg and A. Grothendieck 
it is known that $st_q$ is a finite group for all prime $q \ne p$.
M. Artin obtained the result that the group $st_p$ can be infinite.
Vvedenskii  called this result the Artin effect.\\

\subsubsection{Shafarevich conjecture}
$$H^1 (k,A) \simeq Hom_{ct}(\pi_1(A_k),{\mathbb Q}/{\mathbb Z})  $$
In some cases conjecture proved by Vvedenskii (the cohomologies on the right are taken over continuous cochains). \\

\section{Groups and group scheme}
\subsection{Elements of the theory of algebraic groups and group schemes.}

Let $R$ be a commutative ring with identity. It is known what the affine scheme $Spec \: R$ is \cite{Shafarevich2}.
We recall here, and briefly explain, following \cite{Shafarevich2,Hartshorne}, some concepts related to the class of varieties that are generated by the reduced separated smooth schemes
$(X, {\mathcal O}_X)$ of finite type over an algebraically closed field.
The important notion of separable schema is defined through the concept of the product of schemas and their clousernes.
In turn, the product of schemes is defined as the product of objects
in the category of schemes, but in terms of morphisms of schemes over a basic scheme $S$ (for example, if $S$
is an algebraically closed field) as a fiber product of these morphisms.
A morphism of schemes $\varphi: X \to Y$
 is called a closed embedding if every point $x \in X$ 
 has such
an affine neighborhood $U$ such that the scheme $\varphi^{-1}(U)$
is affine and the homomorphism $\varphi^*: {\mathcal O}_X \to {\mathcal O}_Y$
epimorphic. In the category of schemes over $S$, there is a morphism $(1, 1): X \to X \times_S X$, which is called a diagonal.
A scheme $X$ is called closed if the morphism of its diagonal is a closed embedding, and a scheme over a ring $R$ if the morphism of schemes $R \to Spec \: R$ is given.
A finite group scheme, or a finite group of order $m$ over $R$, is a group scheme locally free of rank $m$ over $R$.
Such a group scheme $G$ is defined by a sheaf of locally free algebras $A$ of  rank $m$ over $R$.
In the works \cite{Shafarevich1}
Serre's quasi-algebraic and pro-algebraic groups are constructed and studied.
In the definition of quasi-algebraic and pro-algebraic groups
according to Serre \cite{Serre1}, the concept of structure is used.
The structure $St$ of a group scheme, or group structure, is given by homomorphisms 1) $\mu: X \times_R X  \to X$  (group law), 2)
$p: X \to X$, $p(x) = e$ (unit), 3) $i: X \to X$ (taking the inverse), satisfying the axioms:
a) $\mu \ast (\mu \times 1) = \mu \ast (1 \times \mu)$ (associativity), b)  $\mu \ast (1, i) =  (i, 1)\ast \mu$ 
 c) $\mu \ast (1, p) =  (p, 1) \ast \mu = 1,$
specifying, respectively, the group law,  taking the inverse element and  unity  which
satisfy the known axioms \cite{Shafarevich2}.
If in an algebraic group $G$, or, more generally, in a group scheme, the group structure $St$ is fixed, then we reveal this through $G_{St}$.

\subsection{Quasi-algebraic and pro-algebraic groups by Serre.}

In what follows, all groups, unless otherwise stated, are assumed to be commutative.
In this section, the letter $K$ denotes a perfect field (algebraically closed field), the letter $p$ denotes its characteristic exponent, that is, $p = 1$ if the characteristic $K$ is equal to zero,
characteristic exponent $= p$ if characteristic $K$ is equal to $p > 0$.
All algebraic varieties are considered defined over $K$. As is known, in the category of algebraic groups
over $K$ there are bijective morphisms that are not morphisms in the sense of algebraic groups. In other words, the category of algebraic groups is additive but not abelian.
Let us recall a well-known example. Let $K = k$ be an algebraically closed field of characteristic
$p > 0$. The topological spaces of algebraic groups $X$ and $Y$, which we denote by the same letters, are given by the condition $X = Y = k$.
Let the group operation of each of the groups be additive and given by the mapping $\mu(x,y) = x + y, x, y \in k$.
Consider the homomorphism $\varphi: X \to Y$ of algebraic groups $X$ and $Y$ given by the condition $\varphi(x) = x^p$.
As a point mapping it is one-to-one and as a mapping of abstract groups it is an isomorphism, but as a regular mapping of manifolds it is not an isomorphism, since   the corresponding ring homomorphism 
$\varphi^{*}: k[Y] \to k[X], \varphi^{*}(T) = T^p,  k[Y] =  k[X] =  k[T], \varphi^{*}( k[Y]) = k[T^p] \neq k[T]$  not is an isomorphism.
\subsubsection{Quasi-algebraic   groups.}
Let $V$ be an algebraic variety and  $\mathcal O$ be the sheaf of functions on $V$.
If $q = p^n, \; n \in {\mathbb Z}$ denote by $\mathcal O^q$  the sheaf whose sections over open sets $U \subset V$
 are the $q$-th powers of the sections of the sheaf $\mathcal O$ over $U$.\\
The concept of a quasi-algebraic group \cite{Serre1} combines into one class algebraic groups between which there are bijections that may not be isomorphisms of algebraic groups.
Let $G_{St}$ be an algebraic group and ${\mathcal O}_{St}$ be a sheaf of functions on $G_{St}$. If $q = p^n, n \in {\mathbb Z}$, then we denote by ${\mathcal O}^q_{St}$ the sheaf whose sections over open sets $U \subset G$
 are the $q$-th powers of the sections of the sheaf ${\mathcal O}_{St}$ over $U$.
 Follow to Serre we also put ${\mathcal O}^{p^{-\infty}} = \bigcup_{n \in {\mathbb Z}}{\mathcal O}^{p^n}.$
 If $q>1$, then the variety $G^q$ corresponding to the sheaf $\mathcal O^q$ is an algebraic group.
 
\begin{proposition}. Let $f: G_1 \to G_2$ be a bijective morphism of algebraic groups. Then there is a positive power
$q$ of $p$, and morphism $g: G_2 \to G_1^q$ such that their composision: 
$$G_1 \stackrel{f}\to G_2 \stackrel{g} \to G_1^q   $$
gives the identity mapping $i: G_1 \to  G_1^q.$  
\end{proposition}
 Let $G$ be a group. If $St$ is an algebraic group structure on $G$ compatible with the group structure, then $G_{St}$ denotes the corresponding algebraic group, and $T_{St}$ and ${\mathcal O}_{St}$ respectively denote the topology and the sheaf of rings.

\begin{proposition}.
\label{prop2}
Let $St_1$ and $St_2$ be two algebraic group structures on $G$ compatible with the group structure.
The following conditions are equivalent:\\
 (i) There is
structure $St_3$ such that the identity mappings $G_{St_3} \to G_{St_1}$ and $G_{St_3} \to G_{St_2}$ are
morphisms.\\
(ii) There is
structure $St_4$ such that the identity mapping $G_{St_4} \to G_{St_1}$ and $G_{St_4} \to G_{St_2}$ are 
morphisms.\\
(iii) There exists for an arbitrary positive power $q$ of the number $p$ the identity mapping $G_{St_1} \to G_{St_2}$ which is a morphism of algebraic groups.\\
(vi) $T_{St_1} = T_{St_2}$ and ${\mathcal O}_{St_1}^{p^{-\infty}} = {\mathcal O}_{St_2}^{p^{-\infty}}$
\end{proposition}

\begin{definition}
\label{def1}
Let $G$ be a group and let  $St_1$ and $St_2$ be two  structeres of algebraic group on $G$ compatible with the group structure. 
Following Serre we say that $St_1$ and $St_2$ are equivalent if they satisfy  Proposition (\ref{prop2}).
 \end{definition}
 Note that, in characteristic zero, equivalence reduces to equality, since $\mathcal O^{p^{-\infty}} = \mathcal O.$
 \begin{definition}
A group $G$ is called a quasi-algebraic group if it is endowed with the class of equivalence (in the sense of definition \ref{def1}) algebraic group structures agreed with its group structure.
 \end{definition} 
 
 \begin{proposition}.
 \label{prop3}
   Let $G$ and $G_1$ be quasi-algebraic groups and let $f: G \to G_1$ be a homomorphism. 
   The following conditions are equivalent:\\
  a)  There are algebraic structures $S$ and $S_1$ on   $G$ and $G_1$ respectively, which are compatible with thair quasi-algebraic structures and such that $f: G \to G_1$ is a morphism of algebraic groups.\\
  b) The mapping $f$  is continuous and if $\varphi$ is a section of $\mathcal O^{p^{-\infty}}_{G_1}$ on an open $U_1$, 
$\varphi \circ f$ is a section of $\mathcal O^{p^{-\infty}}_{G}$  on the open $f^{-1}(U_1)$.\\
  c) The graph of $f$ is a closed subgroup in $G \times G_1$.
 \end{proposition}
   
   \begin{definition}
\label{def3}
Let $G$ and $G_1$ be two quasi-algebraic groups.
A morphism from $G$ to $G_1$ is any homomorphism $f : G \to G_1$ satisfying the equivalent conditions of 
Proposition \ref{prop3}.
 \end{definition}
\begin{proposition}.
 \label{prop4}
Let $f: G \to G_1$ be a morphism of quasi-algebraic groups such that $N$ is the kernel of $f$ and $I$ is its image.
Then $N$ and $I$ are  closed respectively  in $G$ and in $G_1$, and $f$ determines after passing to the factors
an isomorphism of $G/N$ with $I$.
\end{proposition}

Follow to Serre denote by  $\mathcal Z$  the category formed by quasi-algebraic groups and their morphisms.
\begin{proposition}.
 \label{prop5}
The category $\mathcal Z$ is Abelian and the notion of a subobject coincides with the notion of a closed subgroup.
\end{proposition}

\begin{remark}
Recall in connection with Proposition \ref{prop4} axioms that turn an additive category to abelian: \\
AB 1) Every morphism has a kernel and a cokernel.  \\
AB 2)   For every morphism $f$, the canonical morphism from $coim \: f$ to $im\:  f$  is an isomorphism. \\ 
Recall also that   if $char\; k = 0$ then  the category $\mathcal Z$ is identical with the category $\mathcal A$ of algebraic groups.
\end{remark}

\begin{proposition}.
 \label{prop6}
 Every object of the category $\mathcal Z$ is artinian.
\end{proposition}

\begin{example}
 Let $AR$ be an Artinian local ring with algebraically closed field $k$ of characteristic $p > 0$. 
 By M. Atiyah, I. Macdonald, a Noetherian local ring with maximal ideal $\mathfrak m$ is an Artin local ring if 
 ${\mathfrak m}^n = 0$ for some $n$. \\
 Let $K$ be a quasi local field, i.e., a discretely normed field complete with respect
    to the topology induced by the norm and with the algebraically closed residue field $k$, 
    $A$ its ring of valuation, and let $\mathfrak m$ its maximal 
ideal. Let $U = A - \mathfrak m$ be the group of units of $A$ and let $U^n = 1 + {\mathfrak m}^n$.
The ring $A/{\mathfrak m}^n$ is Artinian, whose group of units is identified with the quotient $U/U^{(n)}$. \\
Let $W_n(k)$ be the ring of Witt vectors of length $n$ over $k$;
if $n$ is large enough, we can lift $W_n(k) \to k$ into a homomorphism $W_n(k) \to AR$ that makes $AR$ a $W_n(k)-$module of finite type.
As a module, $AR$ is isomorphic to a direct sum of modules $W_{n_i}(k)$, $n_i \le n$.

As each of the $W_{n_i}(k)$ has a natural structure of algebraic variety on $k$, 
we can transport this structure to $AR$,
and the structure thus obtained does not depend on the choice of the isomophism.\\
By applying the above, we therefore obtain an algebraic group structure on $U/U^{(n)}$,
 and $U$ is the projective limit of the groups $U/U^{(n)}$.
 
 The quotient $U/U^{(1)}$ is identified with the multiplicative group $G_m$. For
$n \ge 1$, the quotient $U^{(n)}/U^{(n+1)}$ is identified with the
additive group $G_a$;

\end{example}

Recall result about  the structure of quasi-algebraic groups.
 \begin{theorem}
Every quasi-algebraic group has a composition series whose successive quotients are isomorphic,
 either to the group $G_a$ or to the group $G_m$, or to an abelian variety, or to a finite group.
\end{theorem}

\begin{remark}
The definition of a quasi-algebraic group can also be given in terms of group schemes.
Let us briefly recall this construction. We extend the category of algebraic groups over $K$ to the category of group schemes over $K$. Since here in what follows we consider only commutative groups, we restrict ourselves to the category of commutative group schemes ${\mathcal CG}_K$ over $K$.
Let $H$ and $G$ lie in  ${\mathcal CG}_K$ and $\varphi: H \to G$ be a purely non-separable isogeny from $H$ to $G$. Let $H$ and $G$ be equivalent if there exists a group scheme $F \in {\mathcal CG}_K$ and purely non-separable isogenies 
$\psi: F \to H, \; \tau: F \to H$. Then a quasi-algebraic group will be a class of equivalent (in the above sense) group schemes.
\end{remark}

\subsubsection{Pro-algebraic   groups.}

Let $G$ be a group with the neutral element $o$, and let $X$ be a homogeneous space on $G$. We will say that $X$ is principal if
 the isotropy subgroup of a point $x \in X$ is reduced to $o$; the choice of $x$ then defines 
a bijection from $G$ onto $X$.\\
  \begin{definition}
Follow to Serre and others we call proalgebraic group a group G endowed with a non-empty family S of subgroups and, for 
all $H \in S$, with a quasi-algebraic group structure on $G/H$, these data satisfying the 
following axioms:\\
$(P_1) \; H, H_1 \in S \Rightarrow H \cap H_1 \in S$.\\
$(P_2)$ If $H \in S$, the subgroups $H_1$ containing $H$ which belong to $S$ are the resiproqed images
 of the closed subgroups of $G/H.$ \\
$(P_3)$  If  $H, H_1 \in S$ and if  $H \subset H_1$, the homomorphism $G/H \to G/H_1$ is the morphism of 
quasi algebraic groups.  \\
$(P_4)$ The natural map $G \to \lim \; inv \; G/H$ (inverse, or projective limit) is a bijection of $G$ onto  the projective limit of groups $G/H \; (H \in S)$.   \\
\end{definition}
 Denote by ${\mathcal PG}_K$ the category of  proalgebraic groups over $K$.
 \begin{definition}
\label{def5}
By Serre and others a group $G \in{\mathcal PG}_K$ is the dimension zero if for any definitive subgroup $H$, the quotient $G/H$ is a finite group.
\end{definition}
\begin{example}
For any prime $p$ the group ${\mathbb Z}_p$ has the dimension zero.
\end{example}

\section{Fundamental groups of  schemes}
\subsection{Homotopy groups.}
In this subsection, we follow to Serre \cite{Serre1} and to  Grothendieck et al. \cite{gav}. Let $G$ be a  quasi-algebraic group. Denote by $G^0$ the connected component of the unit of $G$.  Further $G^0$ is called the connected component of the group $G$. Suppous that $G$ is the  proalgebraic group with complete defining set $\mathcal S$; for $H \in {\mathcal S}$  connected component $(G/H)^0$  of the factor group 
 $G/H$ is closed subgroup in  $G/H$, and, if $H_1 \subset H$ then the image $G/H_1$ in $G/H$ is $(G/H)^0$. In view of this, one can put $G/{G^0} = \lim \; inv \; (G/H)/(G/H)^0$.
 \begin{definition}
\label{def5}
Factor group $G/G^0$ denoted by $\pi_0(G)$ and is called the $0$th homotopy group of the proalgebraic group $G$.
\end{definition}
\begin{remark}
Factorization operation 
$$\pi_0(G) = G/G^0  $$
defines a functor
$$\pi_0: {\mathcal PG}_K \to {\mathcal PG}_K^0  $$
from category ${\mathcal PG}_K$ to category ${\mathcal PG}_K^0$
proalgebraic groups of dimension zero.
\end{remark}
\begin{definition}
\label{def6}
The left derived functors of the functor $\pi_0$ are called the $i$th homotopy groups of the proalgebraic group $G$ and denoted by $\pi_i(G)$.
\end{definition}
\begin{remark}
The presence of a sufficient number of projective objects in the category ${\mathcal PG}_K$ \cite{Serre1} makes Definition \ref{def6} correct.
\end{remark}
\begin{definition}
\label{def7}
 Let  $G \in {\mathcal PG}_K$. The first homotopy group $\pi_1(G)$ of the group $G$ is called the fundamental group of  the group $G$ .
\end{definition}
\begin{definition}
\label{def8}
A group $G\in PG_K$ is called connected if $G = G_0$. Group  $G\in PG_K$ is called singly connected if $\pi_1(G) = 0$.
\end{definition}
\begin{propositiondefinition}
\label{prodef}
Let  $G\in PG_K$. There is the connected and singly connected proalgebraic group $\overline{G}$  and a morphism 
$u: \overline{G} \to G$ such that the kernel and cokernel of $u$ are  proalgebraic groups of dimension zero.
The pair $(\overline{G}, u)$ is unique, up to isomorphism.
The pair $(\overline{G}, u)$ is called the universal covering group of the group $G$.
\end{propositiondefinition}
Recall the example of the   first étale  homotopy group (the  étale fundamental group) $\pi_1(G)$ of the group $G$.
\begin{example}
\label{ea}
 Let  ${\mathbb A}^1(F)$ be the affine line over an algebraically closed field $F$ of characteristic zero and $m$ be the geometric point of  ${\mathbb A}^1(F)$. Then
 $$ \pi_1({\mathbb A}^1 \setminus 0, m) =  \lim \; inv \; \mu_n(F) \simeq {\hat{\mathbb Z}}.$$
\end{example}

\subsection{Fundamental groups of fields}.
In his works, which relate to the arithmetic of number fields and rings, O. N. Vvedenskii used and developed the results of S. Lichtenbaum \cite{Lichtenbaum1}.
In  recent works, C. Lichtenbaum \cite{Lichtenbaum2,Lichtenbaum3} defined the Weil étale topology and the Weil étale site and studied the fundamental groups associated with their.
Let us recall Lichtenbaum's considerations. For  the function field  $K$  of a curve over a finite field $k$ let ${\mathfrak G}_K$ be the Galois group of $\overline K$ over $K$. 
The group $Gal(K{\overline k}/K)$ is isomorphic to $\hat{\mathbb Z}$.
There is a natural surjection $\pi: {\mathfrak G}_K \to Gal(K{\overline k}/K)$.
 The Weil group $W_K \simeq {\pi}^{-1}({\mathbb Z})$ where ${\mathbb Z}$ is the subgroup of $ Gal(K{\overline k}/K$.
 Weil's site and topos are defined in a natural way.\\
 B. Morin's work  \cite{MorinII} follows this circle of ideas.
The studies of these authors, as well as the studies of O. N. Vvedenskii, use the results of A. Grothendieck et al.\cite{gav}.
Since these studies to some extent use and develop the concept of étale topos according to Grothendieck, that is, roughly speaking, the category of étale sheaves on an étale site, we recall the example of an étale sheaf in the situation with (profinite) Galois extensions of finite fields.
\begin{example}
\label{d3}
 Let $X$ be an algebraic variety over ${\mathbb F}_q$ and ${\overline X} = X \times_{{\mathbb F}_q}{\overline {\mathbb F}}_q$. 
  Let ${\overline{\mathfrak g}} = {\mathcal Gal}({\overline {\mathbb F}}_q:{\mathbb F}_q)$ be the corresponding Galois group.
  An etale  sheaf on $X$ corresponds to a sheaf on ${\overline X}$ together  with a continues action of
   ${\overline{\mathfrak g}}$.
\end{example}
In the cited paper \cite{MorinII}, Morin defines the fundamental group underlying the étale Weil cohomology of a number ring.
By Molin,  corresponding Weil étale topos is determined
as a refinemen  of the Weil's site according to Likhtenbaum.
The author of  \cite{MorinII} demonstrates the naturalness of his definition in the case of a smooth projective curve and
further defines the étale fundamental Weil group of an open subscheme of the spectrum of a number ring.
This fundamental group is a projective system of locally compact topological groups that represents the first cohomology with coefficients in a locally compact Abelian group.
This result is usesed in  \cite{MorinII} to calculate cohomology groups of small degrees and to
 check that the Weil étale topos satisfies the expected properties of the Lichtenbaum topos.
 Let $Y$ be an open subscheme of a smooth projective curve over a finite field $k$, and let 
$S_{et}(W_k,{\overline Y})$ be  the topos $W_k-$equivariant etale sheaves on the projective curve
${\overline Y} = Y \otimes_k {\overline k}$.\\

{{\bf Theorem} M1}(\cite{MorinII})  There is an equivalence $Y_W^{sm} \simeq S_{et} (W_k,{\overline Y})$ where $Y_W^{sm}$ is the (small) Weil-\'etale topos defined in this paper.\\

For a connected \'etale ${\overline X}-$scheme ${\overline  U}$ author defines its Weil-\'etale topos as the 
slice topos ${\overline U}_W := {\overline X}_W/{\gamma}^* {\overline U}$. Let $K$ be the number field 
corresponding to the generic point of 
${\overline  U}$, and let $q_{\overline U} : Spec ({\overline K}) \to {\overline U} $ be a geometric point.  
Similarly to the definition of the \'etale fundamental group as a (strict) projective system of finite quotients 
of the Galois group $G_K$, author  defines the analogous (strict) projective system 
$ {\underline W}({\overline U}, q_{\overline U})$ of locally compact quotients of the Weil group $W_K$ .

 {{\bf Theorem} M2}(\cite{MorinII}) The Weil-\'etale topos ${\overline U}_W$ is connected and locally connected over the topos $T$ 
of locally compact spaces. The geometric point $q_{\overline U}$  defines a $T-$valued point $p_{\overline U}$  
of the topos ${\overline U}_W,$ and we have an isomorphism 
$\pi_1 ({\overline U}_W, p_{\overline U } ) \simeq {\underline W}({\overline U}, q_{\overline U}  )$ of 
topological pro-groups.\\

The fundamental group and the fundamental scheme can be defined in terms of the 
 Tannakian category approach (see \cite{Nori,biswassantosi} and references therein).

\section{On Local Fields and  Local Class Field Theory}
 
\subsection{On local fields}
Let $L$ be a finite extension of the local field $K$, $l, k -$ their residue fields, $p = char\; k$,
and $e_{L/K}-$ ramification index of $L$ over $K$.

An extension $L/K$ is called unramified if \\
        a) $e_{L/K}=1$; \\
        b) the extension $ l/k$ is separable. \\
         An extension $e_{L/K}$ is said to be weakly ramified if\\
        a) $p\nmid e_{L/K}$; \\
        b) the extension $l/k$ is separable. \\
An extension $L/K$ is said to be wildly ramified if \\
       $e_{L/K} = [L:K]= (char\; k)^s, s\ge1$.
       
Further, by $Tr_{L/K}$ and $Norm_{L/K}$ we denote, respectively, the trace and the norm
$L/K$ extensions, omitting indices when it is clear which extension is in question.       
       
Denote by $ K_{nr}$ the maximal unramified extension of the field $ K$ (in a fixed algebraic closure of the field $K$) with the residue field $k_s$ , which is the algebraic closure of the field $k$.

Recall that a local field with an algebraically closed residue field is called quasilocal.

\begin{lemma}
\label{pr}
 If local field $K$ contains a primitive root $\xi_p \; p-$th degree of unity, then
$\nu_K (\xi_p - 1) = \frac{e}{p-1}=e_1$ (i.e.) $e_1$ is an integer. Here $e = \nu_K (p), p=char\; k$.
\end{lemma}
 \textbf{Proof.} If $\xi_p$ is a primitive $p-$ root of unity, then ${\xi_p}^p - 1=0$ and $\xi_p$ is a root of a polynomial 
 $x^{(p-1)}+x^{(p-2)} +\ldots +x+1$  irreducible over $ K$ . Then ${\xi_p} - 1$ is the root of the polynomial
$(x+1)^{(p-1)}+(x+1)^{(p-2)}+ \ldots +(x+1)+1 = x^{(p-1)}+ p(\ldots)+p.$
The value of the exponent $p$ at the root of such a polynomial is
$\frac{e}{p-1} $, i.e. $\nu_K (\xi_p - 1) = \frac{e}{p-1}$. The lemma is proven.      
       
\textbf{Corollary.} Under the conditions of Lemma \ref{pr}, $e = \nu_K (p)$ is divisible by $p$.       
       
{\bf Ramification groups.} Let $L/K$ be a finite Galois extension with Galois group
${\mathfrak G} = {\mathfrak G}(L/K)$. Let ${\mathcal O}_K$ be the ring of integers in the field $K$.

We define ramification groups ${\mathfrak G}_i ( i=-1,0,1,\ldots)$ by
${\mathfrak G}_i = \{\sigma \in {\mathfrak G} \mid \nu_L (\sigma a - a)\ge i+1 $ for all
$a \in {\mathcal O}_L \}$.

It is easy to check that the groups ${\mathfrak G}_i$ are normal subgroups of the group
${\mathfrak G}, {\mathfrak G}_{i+1} \subset {\mathfrak G}_i ,{\mathfrak G}_{-1}= {\mathfrak G}$ and
for sufficiently large $m$ ${\mathfrak G}_m = 1$.

Let us introduce a lower and an upper (Herbrand) numbering of ramification groups.
       Let $x$ denote a real variable that is $\ge-1$. Let's put

${\mathfrak G}_x = {\mathfrak G}_i$

 where $ i$ is the smallest integer that is $\ge x$. We introduce the notation

$g_i :=$ (the order of the group ${\mathfrak G}_i$ ).
Let $x$ be a real number and $m$ the integer part of the number $x$.

Define the function
$$\varphi (x) = \left\{
                                   \begin{array}{lcl}(x, if
                                    -1 \le x \le 0 \\
                                   \frac{1}{g_0}  [g_1+ \ldots + g_m +(x-m)g_{m+1} ], if \; x \; \ge 0. \\
                                                                       \end{array}
                                     \right.
                                      $$

The function $\varphi (x)$ is continuous, strictly increasing, and therefore has an inverse function $\psi (y)$, which is also continuous and strictly increasing $(-1 \le y)$. The new, `top' numbering of the ramification groups is now given as follows:

${\mathfrak G}^{\varphi (x)} = {\mathfrak G}_{\psi (y)}$, where $ y= \varphi (x)$ and $x = \psi (y)$ .         \\                             
                                      
{\bf  Different.} Denote by $\pi_K$ the uniformizing element of the field $K$, that is, such
element $\pi_K$ such that $\nu_K (\pi_K )=1$.

Denote by ${\mathfrak m}_K = {\pi}_K \cdot {\mathcal O}_K$ the maximal ideal of the ring ${\mathcal O}_K$. Let $L/K$ be a wildly ramified extension of prime degree $p$. Let us define the  different ${\mathcal D}$ of the extension $L/K$ by the formula
\begin{equation}
\label{dif}
                                                                      {\mathcal D} = (f^{'} ({\pi}_L)),
\end{equation}
where $f(x)$
 is the minimal polynomial for ${\pi}_L$ over $K$. \\
 Note that
 \begin{equation}
\label{difeq}
                {\mathfrak D} \subset {\mathfrak d}\cdot   {\mathcal D}^{-1} \Leftrightarrow Tr({\mathfrak D}) \subset {\mathfrak d},
\end{equation}
where ${\mathfrak d}$ is a fractional ideal in ${\mathcal O}_K$ and ${\mathfrak D}$ is a fractional ideal in ${\mathcal O}_L$.
 
\subsection{On local class field theory}
 Let us illustrate the elements of the local class field theory  and its application on the example of the
invariance of the Hodge-Tate decompositions according to Serre \cite{Serre2}.
Let $K$ be a local field of characteristic zero with perfect residue field $k$ of characteristic $p > O$. 
For   the completion $\mathsf C$  of $\overline K$ the Galois group ${\mathfrak G}_{\overline K}$ acts continuously on $\mathsf C$.\\
In the locally compact case, when  $K$  is the finite extension of ${\mathbb Q}_p$ by the local class field theory 
it is possible to identify ${\mathfrak G}_{K}$ with complation of ${\hat K}^{*}$ of $K^{*}$ and the inertia subgroup of 
${\mathfrak G}_{K}$ with the group of units of $K$.

\section{Vvedenskii`s research}
\subsection{Duality in elliptic curves over a local field}
When $k$ is finite, it is known from Tate results \cite{Tate1} that the group of principal homogeneous
 spaces over $A$ is dual to the group $\overline A$  the Picard variety of the variety $A$ except for the $p$-component, where $p$ is the characteristic of $K$.
 
 When $k$ is algebraically closed, it was shown by Shafarevich \cite{Shafarevich1}  and independently by Ogg \cite{Ogg1} that the group of principal homogeneous spaces over A is dual to the group $\pi_1({\overline A}_K)$, i.e., the fundamental group of the
 pro-algebraic group ${\overline A}_K$ in the sense by Serre \cite{Serre1} except for the $p$-component, where $p$ is the characteristic of $k$.
 
 It was conjectured that the duality holds also 
for the $p$-component. In the present work \cite{Vvedenskii64} this conjecture is proved for the special case when $A$ is an elliptic curve  whose reduction has a Hasse invariant other than $0$.

Vvedenskii remarks that he do not find explicitly a natural pairing between $\pi_1({\overline A}_K)$ and the group of principal homogeneous spaces, although the proof of the duality, which  is done in a purely computational way, permits him to deduce that one exists.\\

\end{document}